\documentclass{article}             
\usepackage{amssymb}
\usepackage{graphicx}
\def\bR{{\mathbb{R}}}
\def\fg{{\mathfrak{g}}}
\def\fh{{\mathfrak{h}}}
\def\fm{{\mathfrak{m}}}

\def\rad{{\mathrm{rad}}}

\newtheorem       {theorem}{Theorem}

\newtheorem{lemma}[theorem]{Lemma}

\begin{document}

\begin{center}
{\Large{The minimal number of homogeneous geodesics\\
depending on the signature of the Killing form}}
\bigskip

{\large{Zden\v ek Du\v sek}}
\bigskip

Institute of Technology and Business in \v Cesk\'e Bud\v ejovice\\
Okru\v zn\'\i\ 517/10, 370 01 \v Cesk\'e Bud\v ejovice, Czech Republic\\
zdusek@mail.vstecb.cz
\end {center}
\bigskip

\begin{abstract} 
The existence of at least two homogeneous geodesics in any homogeneous Finsler manifold
was proved in a previous paper by the author.
The examples of solvable Lie groups with invariant Finsler metric which admit just two homogeneous
geodesics were presented in another paper.
In the present work, it is shown that a homogeneous Finsler manifold with indefinite Killing form
admits at least four homogeneous geodesics.
Examples of invariant Randers metrics on Lie groups with definite Killing form
admitting just two homogeneous geodesics
and examples with indefinite Killing form admitting just four homogeneous geodesics are presented.
\end{abstract}
\bigskip

\noindent
{\bf MSClassification:} {53C22, 53C60, 53C30}\\
{\bf Keywords:} {Homogeneous Finsler manifold, homogeneous geodesic}

\section{Introduction}
The existence of at least one homogeneous geodesics in arbitrary homogeneous Riemannian manifold
was proved by O. Kowalski and J. Szenthe in \cite{KS}, by an algebraic construction in the Lie algebra.
In the papers \cite{KNV} and \cite{KVl}, it was proved that this result is optimal, namely,
examples of homogeneous Riemannian metrics on solvable Lie groups were constructed
which admit just one homogeneous geodesic through any point.
Generalization of this existence result to pseuddo-Riemannian geometry was proved by the author
in \cite{Djgp}, in the more general framework of affine geometry, using a purely affine approach and differential topology.

Generalization of the above existence result to Finsler geometry was proved in the series of papers
\cite{YD1} by Z. Yan and S. Deng for Randers metrics (by the algebraic construction),
\cite{Due1} by the author for odd-dimensional Finsler metrics,
\cite{Due2} by the author for Berwald or reversible Finsler metrics
(in both cases using the affine approach),
\cite{YH} by Z. Yan and L. Huang in general
(using again the original idea by Kowalski and Szenthe and a purely Finslerian construction).
However, due to the nonreversibility of general Finsler metrics, it was conjectured by the author
that the result and its proofs in the nonreversible situation are not optimal, namely that an arbitrary
homogeneous Finsler manifold admits at least two homogeneous geodesics through arbitrary point.

In comparison with Riemannian geometry, the situation is rather delicate.
In the context of Finsler geometry, the trajectory of the unique homogeneous geodesic in a Riemannian manifold
should be regarded as two geodesics, with initial vectors $X$ and $-X$,
and only reparametrizations in the same direction may be identified as a trajectory of a geodesic.
For a general homogeneous Finsler manifold, the initial vectors of the two homogeneous geodesics may be non-opposite.
In \cite{DuSr2}, examples of invariant Randers metrics which admit just two homogeneous geodesics
are constructed. The initial vectors of these geodesics are $X+Y$ and $-X+Y$, for certain vectors $X,Y\in T_pM$.
For the construction, Randers metrics which are modifications of Riemannian metrics
of examples from \cite{KNV} and \cite{KVl} are used.
These examples are the solvable Lie groups and hence their Killing form vanishes identically.
It was also demonstrated with an example that general Randers metrics whose underlying Riemannian metric
admits just two homogeneous geodesics (with initial vectors $X$ and $-X$) may admit more than two homogeneous geodesics.

The proof of the existence of at least two homogeneous geodesics in Finsler geometry was given in \cite{DuS},
by the geometrical interpretation of the algebraic condition from the original proof by Kowalski and Szenthe
and used later by Yan and Huang.
The proof in \cite{DuS} contained a small inaccuracy, namely, it was implicitly assumed that the Killing form
of the isometry group is negatively semidefinite.
If the restriction of the Killing form to $\fm$ is indefinite, the proof works as well.
Moreover, there are at least two further solutions in this situation.
In the present paper, in Section 4, we give a full proof of the existence of four homogeneous geodesics
for arbitrary homogeneous Finsler manifold such that the restriction of the Killing form to $\fm$ is indefinite.
Before this, in Section 3, we illustrate the crucial geometrical idea with the examples of Randers metrics
on Lie groups in dimension $3$ and with definite Killing form admitting just two homogeneous geodesics
and Randers metrics with indefinite Killing form admitting just four homogeneous geodesics.

\section{Preliminaries}
Let $(M,F)$ be a Finsler manifold.
If there is a connected Lie group $G$ which acts transitively on $M$ as a group of isometries,
then $M$ is called a {\it homogeneous manifold\/}.
Homogeneous manifold $M$ can be naturally identified
with the {\it homogeneous space\/} $G/H$, where $H$ is the isotropy group of the origin $p\in M$.
A homogeneous Finsler space $(G/H,F)$ is always a {\it reductive homogeneous space\/}:
We denote by $\fg$ and $\fh$ the Lie algebras of $G$ and $H$ respectively
and consider the adjoint representation ${\mathrm{Ad}}\colon H\times\fg\rightarrow\fg$ of $H$ on $\fg$.
There exists a {\it reductive decomposition} of the form $\fg=\fm+\fh$ where $\fm\subset\fg$
is a vector subspace such that ${\rm Ad}(H)(\fm)\subset\fm$.
For a fixed reductive decomposition $\fg=\fm+\fh$ there is the natural identification
of $\fm\subset\fg=T_eG$ with the tangent space $T_pM$ via the projection $\pi\colon G\rightarrow G/H=M$.
Using this natural identification, from the Minkovski norm and its fundamental tensor on $T_pM$,
we obtain the ${\mathrm{Ad}}(H)$-invariant Minkowski norm and the ${\mathrm{Ad}}(H)$-invariant fundamental tensor
on $\fm$ and we denote these again by $F$ and $g$.

Special Minkowski norms (on a vector space ${\mathbb{V}}$) are the {\it Randers norms}.
They are determined by a symmetric positively definite bilinear form $\alpha$
and a vector $V$ such that $\alpha(V,V)<1$, or, equivalently, its $\alpha$-equivalent 1-form $\beta$
related with $V$ by the formula
$\beta(U) = \alpha(V,U)$ for all $U\in {\mathbb{V}}$.
The Randers norm $F$ is then defined by the formula
\begin{eqnarray}
\label{randers}
F(U) = \sqrt{\alpha(U,U)} + \beta(U),\qquad \forall U\in {\mathbb{V}}.
\end{eqnarray}
If a Finsler metric $F$ on $M$ restricted to any tangent space $T_pM$ is a Randers norm, it is called a {\it Randers metric}.
Obviously, a Randers metric $F$ is determined by a Riemannian metric $\alpha$ and a smooth $1$-form $\beta$.
A homogeneous Randers metric $F$ is determined by a Randers norm on $\fm$,
in other words by a symmetric positively definite $2$-form and a $1$-form on $\fm$
and these forms are denoted again by $\alpha$ and $\beta$.
We remark that, in the literature, the letter $\alpha$ is sometimes used for the norm induced by the $2$-form $\alpha$
and then formula (\ref{randers}) above is without the square root.
We choose the notation above because for $\beta=0$, $F$ is the Riemannian norm and components $g_{ij}$
of the fundamental tensor are just the components of the Riemannian metric $\alpha$.

We further recall that the {\it slit tangent bundle} $TM_0$ is defined as $TM_0=TM\setminus\{ 0\}$.
Using the restriction of the natural projection $\pi\colon TM\rightarrow M$ to $TM_0$,
we naturally construct the pullback vector bundle $\pi^*TM$ over $TM_0$.
The {\it Chern connection} is the unique linear connection on the vector bundle
$\pi^*TM$ which is torsion free and almost $g$-compatible,
see some monograph, for example \cite{BCS} by D. Bao, S.-S. Chern and Z. Shen or \cite{De} by S. Deng for details.
Using the Chern connection, the derivative along a curve $\gamma(t)$ can be defined.
A regular smooth curve $\gamma$ with tangent vector field $T$ is a {\it geodesic} if $D_{T} (\frac{T}{F(T)}) = 0$.
In particular, a geodesic of constant speed satisfies $D_{T} {T} = 0$.
A geodesic $\gamma(s)$ through the point $p$ is {\it homogeneous} if it is an orbit
of a one-parameter group of isometries.
More explicitly, if there exists a nonzero vector $X\in\fg$ such that $\gamma(t)={\rm exp}(tX)(p)$ for all $t\in\bR$.
The vector $X$ is called a {\it geodesic vector}.
Geodesic vectors are characterized by the following {\it geodesic lemma}.
\begin{lemma}[\cite{La}]
\label{golema2}
Let $(G/H,F)$ be a homogeneous Finsler space with a reductive decomposition $\fg=\fm+\fh$.
A nonzero vector $Y\in{\fg}$ is geodesic vector if and only if it holds
\begin{equation}
\label{gl2}
g_{Y_\fm} ( Y_{\mathfrak m}, [Y,U]_{\mathfrak m} ) = 0 \qquad \forall U\in{\mathfrak m},
\end{equation}
where the subscript $\fm$ indicates the projection of a vector from $\fg$ to $\fm$.
\end{lemma}
We shall use this lemma for Randers metrics $F=\sqrt{\alpha}+\beta$
and in the situation with trivial algebra $\fh$.
In this special situation, the above statement has the following form.
\begin{lemma}[\cite{DuC}]
\label{c1}
Let $F=\sqrt{\alpha}+\beta$ be a homogeneous Randers metric on $G$, let $\fg$ be the Lie algebra of $G$
and $V\in\fm$ be the vector $\alpha$-equivalent with $\beta$.
The vector $X\in \fg$ is geodesic if and only if
\begin{eqnarray}
\label{golemanew}
\alpha \Bigl ( X + {\sqrt{\alpha(X,X)}}\cdot  V, [X,U] \Bigr ) = 0 \qquad \forall U\in\fm.
\end{eqnarray}
\end{lemma}

In the proof of the existence of a homogeneous geodesic, see \cite{KS}, \cite{YH} or \cite{DuS},
the geometrical interpretation of the crucial idea is the following.
With respect to the Killing form $K$, the bracket $[Y,U]_\fm$ is always orthogonal to $Y_\fm$.
We shall denote the restriction of $K$ from $\fg$ to $\fm$ again by $K$
and we shall consider orthogonal complements in the vector space $\fm$.
If we find a vector $Y\in\fg$ such that $K(Y_\fm)\neq 0$ and the orthogonal spaces to $Y_\fm$ with respect to $g_{Y_\fm}$
and with respect to $K$ are equal, then the bracket $[Y,U]_\fm$ will be also orthogonal to $Y_\fm$
with respect to $g_{Y_\fm}$, we reach the equality
\begin{eqnarray}
\label{c4}
g_{Y_\fm} ( Y_{\mathfrak m}, [Y,U]_{\mathfrak m} ) = 
K         ( Y_{\mathfrak m}, [Y,U]_{\mathfrak m} ) = 0 \qquad \forall U\in{\mathfrak m}
\end{eqnarray}
and, according to geodesic lemma, $Y$ will be geodesic vector.
We denote by $I_F$ the unit indicatrix in $\fm$ given by the condition $F(X)=1$ and we denote by
$S_K$ the unit (pseudo-)sphere in $\fm$ given by the condition $K(X,X)=\pm 1$.
We shall use the geometrical property that the orthogonal space to a vector $X\in I_F$
with respect to the scalar product $g_X$ is the tangent space to the indicatrix $I_F$ at $X$.
In the same way, the orthogonal space to a vector $X\in S_K$
with respect to $K$ is the tangent space to the (pseudo-)sphere $S_K$ at $X$.
For each vector $X\in\fm$ such that $K(X,X)\neq 0$, we put $X_F=X/F(X)\in I_F$ and $X_K=X/| K(X,X)| \in S_K$.
We are looking for vectors $X\in\fm$, such that the tangent space to $I_F$ at $X_F$
and the tangent space to $S_K$ at $X_K$ are equal.
If we interpret these tangent spaces (vector subspaces of $\fm$) as affine subspaces of $\fm$
and put them into their origin points ($X_F$, or $X_K$, respectively), these spaces will be parallel.
Such vectors $X\in\fm$ will satisfy the above conditon (\ref{c4}).

We shall illustrate the situation with examples of Lie groups with invariant Randers metrics
and whose $\rad(K)$ is trivial.
In the first example, $G={\mathrm{SO}}(3)$, the Killing form is definite and $G$ admits just two homogeneous geodesics,
in generic situation.
In the second example, $G={\mathrm{SL}}(2)$, the Killing form is indefinite and $G$ admits just four homogeneous geodesics,
in generic situation.
Further, we shall generalize the procedure of finding at least two homogeneous geodesics
in any homogeneous Finsler manifold described in \cite{DuS}.
We show that any homogeneous Finsler manifold
such that the restriction of the Killing form of the isometry group $G$ from $\fg$ to $\fm$
is indefinite admits at least four homogeneous geodesics.

\section{Examples}
\subsection{Example 1, $\fg_1\simeq{\mathfrak{so}}(3)$}
Consider the Lie algebra $\fg_1={\mathrm{span}}\{E_i\}_{i=1}^3$ generated by the Lie brackets
\begin{eqnarray}
\nonumber
[E_1,E_2]  =   a E_3, \qquad
[E_1,E_3]  =  -b E_2, \qquad
[E_2,E_3]  =   c E_1.
\end{eqnarray}
In the matrix form, for the special choice $a=b=c=1$, we can identify the generators $E_i$ with the matrices
\begin{eqnarray*}
E_1= \left ( \begin{array}{ccc}
 0 & 1 & 0 \cr
-1 & 0 & 0 \cr
 0 & 0 & 0
\end{array} \right ), \quad
E_2= \left ( \begin{array}{ccc}
 0 & 0 & 0 \cr
 0 & 0 & 1 \cr
 0 &-1 & 0
\end{array} \right ), \quad
E_3= \left ( \begin{array}{ccc}
 0 & 0 & 1 \cr
 0 & 0 & 0 \cr
-1 & 0 & 0
\end{array} \right ).
\end{eqnarray*}
It is easy to check that $\fg_1\simeq{\mathfrak{so}}(3)$ for $a,b,c>0$.
By the direct calculations, we also easily check that the Killing form, with respect to the basis $B=\{E_1,E_2,E_3\}$, is
\begin{eqnarray}
\label{KF}
K = -2ab\, x_1^2 -2ac\, x_2^2 -2bc\, x_3^2.
\end{eqnarray}
We now put $X=x_1E_1+x_2E_2+x_3E_3$ and we write down the Lie brackets
\begin{eqnarray}
\nonumber
 [ X,E_{1} ] & = &   - a x_2 E_3 + b x_3 E_2 , \cr
 [ X,E_{2} ] & = &     a x_1 E_3 - c x_3 E_1 , \cr
 [ X,E_{3} ] & = &   - b x_1 E_2 + c x_2 E_1.
\end{eqnarray}
From Lemma \ref{c1} and the equation (\ref{golemanew}) with
$\alpha$ given by the identity matrix with respect to the basis $B$ above and
with $V=v_1E_1+v_2E_2+v_3E_3$, we obtain the system of equations
\begin{eqnarray}
\nonumber
 bx_3 \left (x_2 + {\sqrt{\alpha(X,X)} v_2} \right ) - ax_2 \left (x_3 + {\sqrt{\alpha(X,X)} v_3} \right ) & = & 0, \cr
-cx_3 \left (x_1 + {\sqrt{\alpha(X,X)} v_1} \right ) + ax_1 \left (x_3 + {\sqrt{\alpha(X,X)} v_3} \right ) & = & 0, \cr
 cx_2 \left (x_1 + {\sqrt{\alpha(X,X)} v_1} \right ) - bx_1 \left (x_2 + {\sqrt{\alpha(X,X)} v_2} \right ) & = & 0,
\end{eqnarray}
which simplifies into the form
\begin{eqnarray}
\label{sysabc0}
\left ( b - a \right ) x_2 x_3 + \left ( b x_3 v_2 - a x_2 v_3 \right ) {\sqrt{\alpha(X,X)}} & = & 0, \cr
\left ( a - c \right ) x_1 x_3 + \left ( a x_1 v_3 - c x_3 v_1 \right ) {\sqrt{\alpha(X,X)}} & = & 0, \cr
\left ( c - b \right ) x_1 x_2 + \left ( c x_2 v_1 - b x_1 v_2 \right ) {\sqrt{\alpha(X,X)}} & = & 0.
\end{eqnarray}

Let us first investigate the situation with Riemannian metrics determined by the orthonormal basis $B$, which means $v_i=0$.
The system of equations simplifies further into the form
\begin{eqnarray}
\label{sysabc}
\left ( b - a \right ) x_2 x_3  & = & 0, \cr
\left ( a - c \right ) x_1 x_3  & = & 0, \cr
\left ( c - b \right ) x_1 x_2  & = & 0.
\end{eqnarray}
For $a=b=c$, any vector $X=(x_1,x_2,x_3)$ is obviously a solution of this system.
From the geometrical point of view, we observe that the unit indicatrix $I_F$
given by the (Riemannian) Finsler function $F = \sqrt{x_1^2+x_2^2}$ is just the coordinate sphere with radius $1$.
The hypersurface $S_K$ 
is a coordinate sphere with radius $\frac{\sqrt{2}}{2a}$.
Hence, obviously, for each $0\neq X\in\fm$, the tangent plane to $I_F$ at $X_F$ and the tangent plane
to $S_K$ at $X_K$ are parallel.

For different values of parameters $a,b,c$, the unit indicatrix $I_F$ is still the coordinate unit sphere
(green in the picture which follows),
but the hypersurface $S_K$ is the coordinate elipsoid (blue in the picture).
Each nonzero vector $X\in\fm$ determines vectors $X_F\in I_F$ and $X_K\in S_K$.
The vectors $X$ for which the tangent space to the indicatrix $I_F$ at $X_F$ is parallel with the tangent space
to the hypersurface $S_K$ at $X_K$ are the positive multiples of coordinate vectors (red in the picture).
If two of the parameters $a,b,c$ are equal, then also positive multiples of all vectors
in the corresponding coordinate plane have this property.
We illustrate the situation with the picture in the coordinate plane $x_3=0$ with $a=b=\frac{1}{2}$, $c=3$.
For the Finsler function $F = \sqrt{x_1^2+x_2^2}$, the indicatrix $I_F$ is just the coordinate unit sphere,
the Killing form is
\begin{eqnarray}
K = -\frac{1}{2}\, x_1^2 - 3\, x_2^2
\end{eqnarray}
and the hypersurface $S_K$ is the blue ellipse in the picture.
The nontrivial solutions of the system (\ref{sysabc}), with the restriction $x_3=0$,
are just the positive multiples of vectors
$X_1 = (1,0)$,
$X_2 = (-1,0)$,
$X_3 = (0,1)$,
$X_4 = (0,-1)$.

\hspace{1 cm}
\includegraphics[width=8cm]{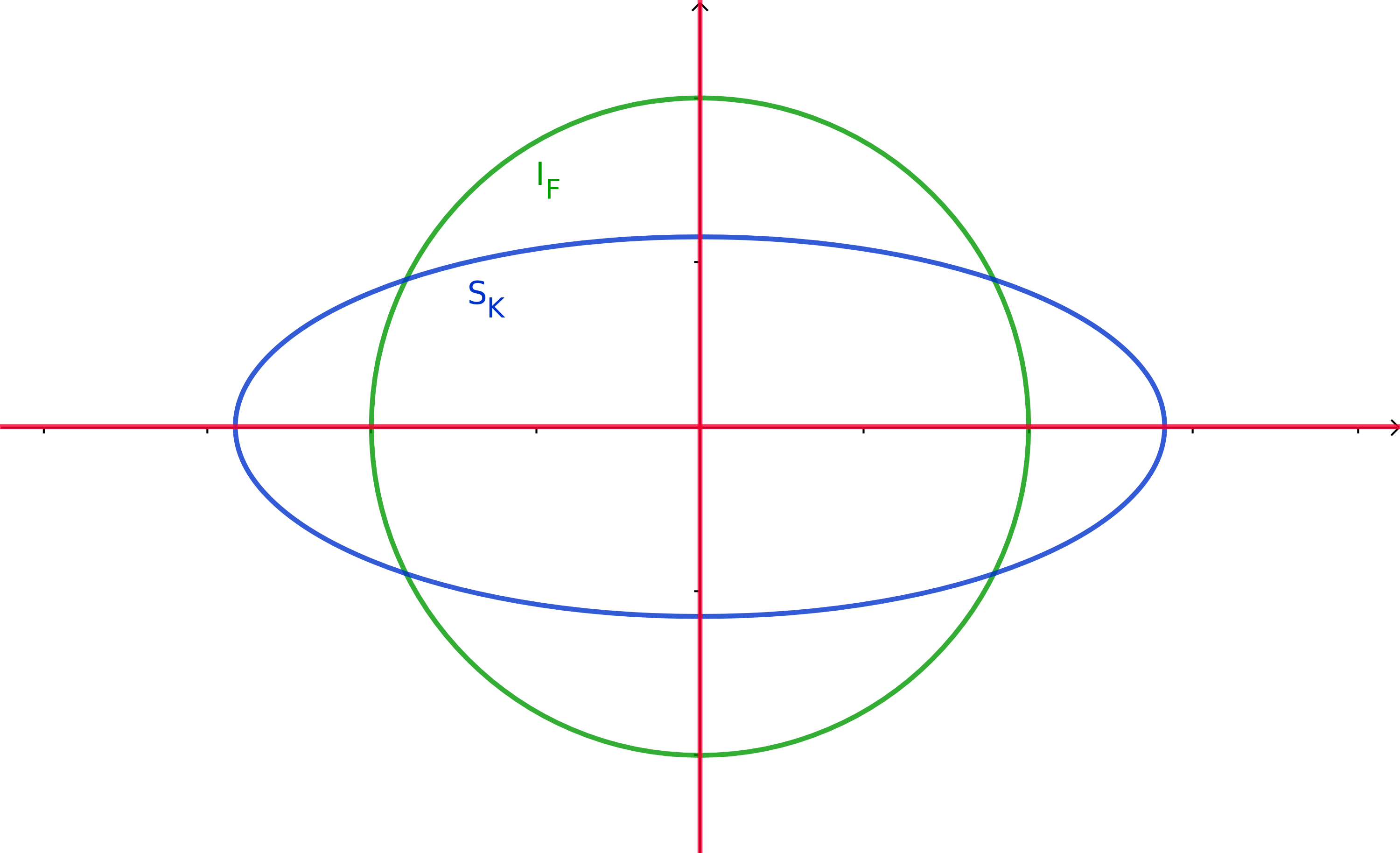}

Let us now turn to general Randers metrics.
For the simplicity, let us consider $V=v_1E_1$, hence $v_2=v_3=0$.
The system of equations (\ref{sysabc0}) simplifies into the form
\begin{eqnarray}
\nonumber
\left ( b - a \right ) x_2 x_3                                   & = & 0, \cr
\Bigl ( \left ( a - c \right ) x_1 - c v_1  {\sqrt{\alpha(X,X)}} \Bigr ) x_3 & = & 0, \cr
\Bigl ( \left ( c - b \right ) x_1 + c v_1  {\sqrt{\alpha(X,X)}} \Bigr ) x_2 & = & 0.
\end{eqnarray}
If $a\neq b$, the solutions are:\\
1) $x_2=x_3=0$,\\
2) $x_2=0$, $x_3\neq 0$ and possible solutions of the equation
\begin{eqnarray}
\label{eqF1}
( a - c ) x_1 = {c v_1} {\sqrt{\alpha(X,X)}},
\end{eqnarray}
3) $x_3=0$, $x_2\neq 0$ and possible solutions of the equation
\begin{eqnarray}
\label{eqF2}
( b-c ) x_1 = {c v_1} {\sqrt{\alpha(X,X)}}.
\end{eqnarray}
Let us illustrate the situation for the particular values of the parameters,
which we set for example $v_1=\frac{1}{2}$, $a=c=1$, $b=2$.
The equation (\ref{eqF1}) has no solutions.
We look for solutions of equation (\ref{eqF2}) such that $\alpha(X,X)=1$ and we obtain
\begin{eqnarray} \nonumber
x_1 = \frac{c v_1}{ b-c } = \frac{1}{2}, \qquad x_2 = \pm\frac{\sqrt{3}}{2}.
\end{eqnarray}
Altogether, we have the four solutions, up to positive scalar multiples.
With respect to the basis $B$, they are
\begin{eqnarray} \nonumber
X_1 & = & (1,0,0),\cr
X_2 & = & (-1,0,0),\cr
X_3 & = & (\frac{1}{2},\frac{\sqrt{3}}{2},0),\cr
X_3 & = & (\frac{1}{2},-\frac{\sqrt{3}}{2},0).
\end{eqnarray}
Again, we illustrate these solutions in the plane $x_3=0$
with a picture. The Finsler function (with the restriction $x_3=0$) is
\begin{eqnarray}
\label{specF}
F = \sqrt{x_1^2+x_2^2} + \frac{1}{2}x_1
\end{eqnarray}
and the hypersurface $I_F$ is given by the equation $F(X)=1$.
The Killing form is given by the formula (\ref{KF}).
For our values of parameters and in the plane $x_3=0$, the hypersurface $S_K$ is given
by the equation
\begin{eqnarray}
2 x_1^2 + x_2^2 = 1.
\end{eqnarray}
In the picture, directions given by the positive multiples of solutions $X_i$ above are in red.
We can see geometrically that the tangent space to $I_F$ at ${X_i^F}$
and the tangent space to $S_K$ at ${X_i^K}$ are parallel, for each $i=1,\dots,4$.

\hspace{1 cm}
\includegraphics[width=8cm]{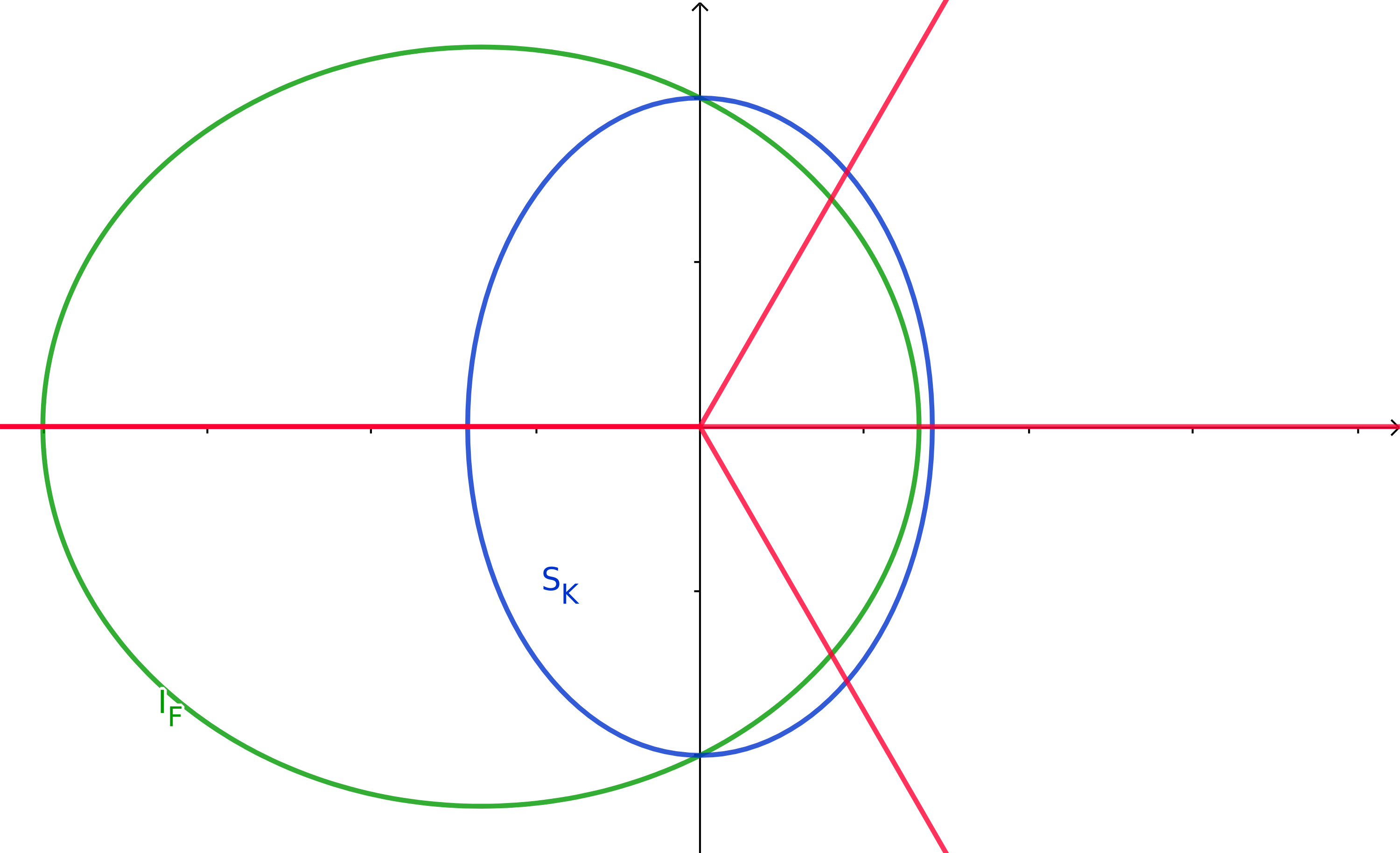}

An important observation shows that the situation may be different for another values of the parameters.
For example, let us choose $v_1=\frac{1}{2}, a=b=2, c=3$.
Because $|x_1| < \alpha(X,X)$ and ${c v_1} > 1$,
none of the equations (\ref{eqF1}) and (\ref{eqF2}) have any nonzero solution
and we are left with solutions $X_1=(1,0,0)$ and $X_2=(-1,0,0)$.
Again, we illustrate the situation with a picture in the plane $x_3=0$.
The hypersufrace $I_F$ is the same as before, the hypersufrace $S_K$ is given by the equation
\begin{eqnarray}
4 x_1^2 + 6 x_2^2 = 1.
\end{eqnarray}
We see also geometrically from the picture that there are no nontrivial solutions $X$
(other than those in the direction of the $x_1$-axis)
in this plane such that the tangent space to $I_F$ at ${X^F}$ and the tangent space to $S_K$ at ${X^K}$ are parallel.

\hspace{1 cm}
\includegraphics[width=8cm]{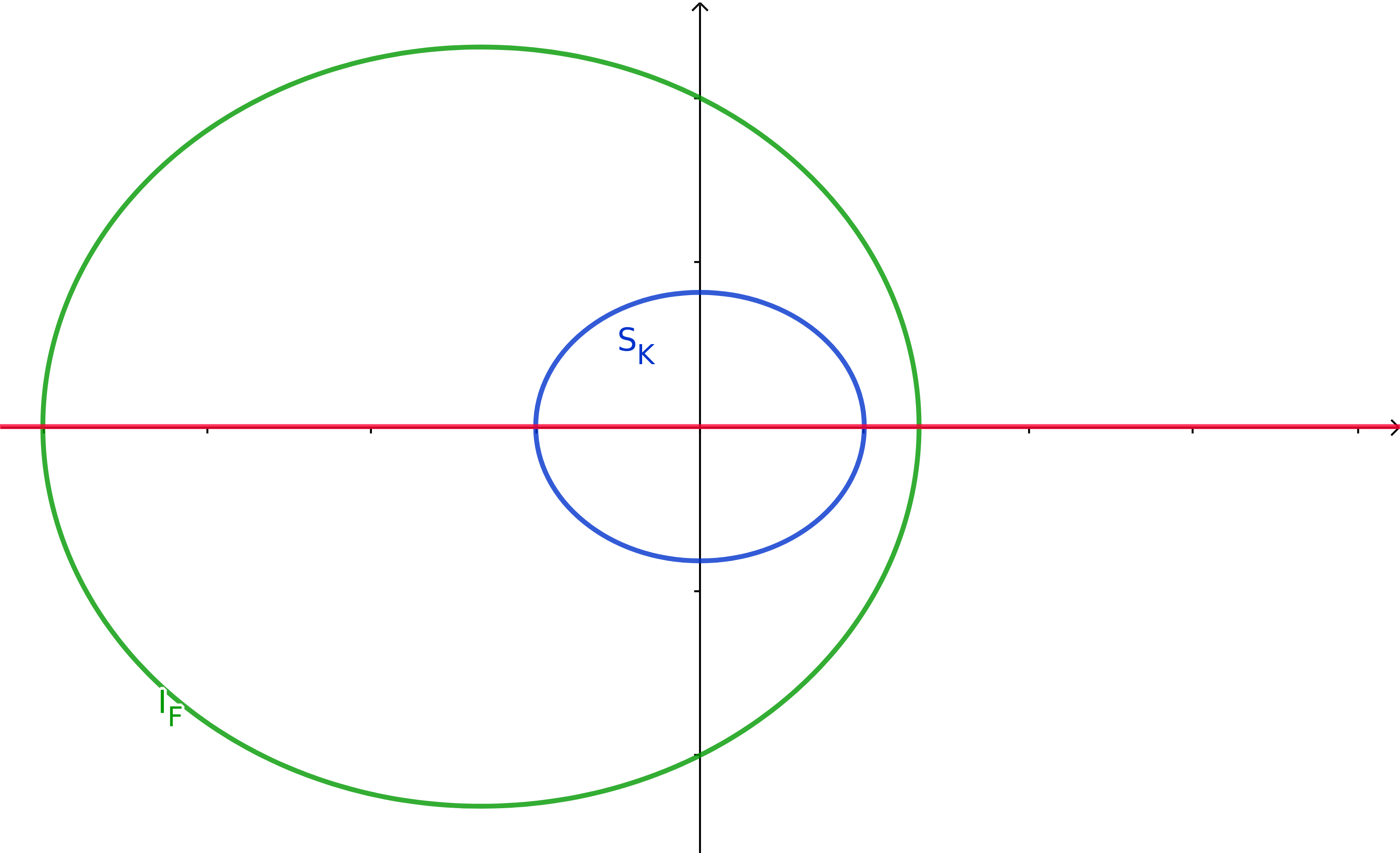}

\subsection{Example 2, $\fg_2\simeq{\mathfrak{sl}}(2)$}
Consider the Lie algebra $\fg_2$ generated by the Lie brackets
\begin{eqnarray}
\nonumber
[E_1,E_2]  =   a E_3, \qquad
[E_1,E_3]  =   b E_2, \qquad
[E_2,E_3]  =   c E_1.
\end{eqnarray}
In the matrix form, for the special choice $a=b=c=1$, we can identify the generators $E_i$ with the matrices
\begin{eqnarray*}
E_1= \left ( \begin{array}{cc}
 0 & 1 \cr
-1 & 0
\end{array} \right ), \quad
E_2= \left ( \begin{array}{cc}
 0 & 1 \cr
 1 & 0
\end{array} \right ), \quad
E_3= \left ( \begin{array}{cc}
 1 & 0 \cr
 0 &-1
\end{array} \right ).
\end{eqnarray*}
It is easy to check that $\fg_2\simeq{\mathfrak{sl}}(2)$ for $a,b,c>0$.
By the direct calculations we obtain that the Killing form,
with respect to the basis $B=\{E_1,E_2,E_3\}$, is
\begin{eqnarray}
K = 2ab\, x_1^2 -2ac\, x_2^2 + 2bc\, x_3^2.
\end{eqnarray}
We put again $X=x_1E_1+x_2E_2+x_3E_3$ and we write down the Lie brackets
\begin{eqnarray}
\nonumber
 [ X,E_{1} ] & = &   - a x_2 E_3 - b x_3 E_2 , \cr
 [ X,E_{2} ] & = &     a x_1 E_3 - c x_3 E_1 , \cr
 [ X,E_{3} ] & = &     b x_1 E_2 + c x_2 E_1.
\end{eqnarray}
From Lemma \ref{c1} and the equation (\ref{golemanew}) with
$\alpha$ given by the identity matrix with respect to the basis $B$ above and
with $V=v_1E_1+v_2E_2+v_3E_3$, we obtain the~system of equations
\begin{eqnarray}
\nonumber
 bx_3 \left (x_2 + {\sqrt{\alpha(X,X)} v_2} \right ) + ax_2 \left (x_3 + {\sqrt{\alpha(X,X)} v_3} \right ) & = & 0, \cr
 cx_3 \left (x_1 + {\sqrt{\alpha(X,X)} v_1} \right ) - ax_1 \left (x_3 + {\sqrt{\alpha(X,X)} v_3} \right ) & = & 0, \cr
 cx_2 \left (x_1 + {\sqrt{\alpha(X,X)} v_1} \right ) + bx_1 \left (x_2 + {\sqrt{\alpha(X,X)} v_2} \right ) & = & 0,
\end{eqnarray}
which simplifies into the form
\begin{eqnarray}
\nonumber
\left ( a + b \right ) x_2 x_3 + \left ( b x_3 v_2 + a x_2 v_3 \right ) {\sqrt{\alpha(X,X)}} & = & 0, \cr
\left ( a - c \right ) x_1 x_3 + \left ( a x_1 v_3 - c x_3 v_1 \right ) {\sqrt{\alpha(X,X)}} & = & 0, \cr
\left ( b + c \right ) x_1 x_2 + \left ( c x_2 v_1 + b x_1 v_2 \right ) {\sqrt{\alpha(X,X)}} & = & 0.
\end{eqnarray}
For the simplicity, let us consider again just the case $V=v_1E_1$, hence $v_2=v_3=0$.
This system of equations simplifies further into the form
\begin{eqnarray}
\nonumber
\left ( a + b \right ) x_2 x_3             & = & 0, \cr
\Bigl ( \left ( a - c \right ) x_1 - c v_1 {\sqrt{\alpha(X,X)}} \Bigr ) x_3 & = & 0, \cr
\Bigl ( \left ( b + c \right ) x_1 + c v_1 {\sqrt{\alpha(X,X)}} \Bigr ) x_2 & = & 0.
\end{eqnarray}
The solutions are:\\
1) $x_2=x_3=0$,\\
2) $x_2=0$, $x_3\neq 0$ and solutions of the equation
\begin{eqnarray}
\label{eqF3}
( a-c ) x_1 = c v_1 {\sqrt{\alpha(X,X)}},
\end{eqnarray}
3) $x_3=0$, $x_2\neq 0$ and
\begin{eqnarray}
\label{eqF4}
x_1 = \frac{ - c v_1}{ b + c } {\sqrt{\alpha(X,X)}}.
\end{eqnarray}
Let us notice that 
$ |v_1| < 1, |\frac{c}{b+c}| < 1$ and hence, with the assumption $\alpha(X,X)=1$, there is always
a solution of the equation (\ref{eqF4}) with $|x_1| < 1$ and $x_2$ determined from the assumption $\alpha(X,X)=1$.
On the other hand, the equation (\ref{eqF3}) may have no nontrivial solutions.
Let us set the particular values of the parameters, for example $v_1=\frac{1}{2}$, $a=b=c=1$.
The equation (\ref{eqF3}) has no nontrivial solution and
the solution of the equation (\ref{eqF4}) is (up to a positive multiple)
\begin{eqnarray} \nonumber
x_1 = \frac{- c v_1}{ b+c } = -\frac{1}{4}, \qquad x_2 = \pm\frac{\sqrt{15}}{4}.
\end{eqnarray}
Altogether, we have the four solutions, up to positive scalar multiples.
With respect to the basis $B$, they are
\begin{eqnarray} \nonumber
X_1 & = & (1,0,0),\cr
X_2 & = & (-1,0,0),\cr
X_3 & = & (-\frac{1}{4},\frac{\sqrt{15}}{4},0),\cr
X_3 & = & (-\frac{1}{4},-\frac{\sqrt{15}}{4},0).
\end{eqnarray}
Again, we illustrate these solutions, which are all in the plane $x_3=0$, with a picture.
The Finsler function (with the restriction $x_3=0$) is still given by the equation (\ref{specF})
and the hypersurface $I_F$ is given by the equation $F(X)=1$.
For our values of the parameter and in the plane $x_3=0$, the hypersurface $S_K$ is given
by the equation
\begin{eqnarray}
2 x_1^2 - 2 x_2^2 = 1.
\end{eqnarray}
In the picture, directions given by the positive multiples of solutions $X_i$ above are again in red.
Again, we can see that, for each $i=1,\dots,4$,
the tangent space to $I_F$ at ${X_i^F}$ and the tangent space to $S_K$ at ${X_i^K}$ are parallel.

\vspace{.4 cm}
\hspace{1 cm}
\includegraphics[width=8cm]{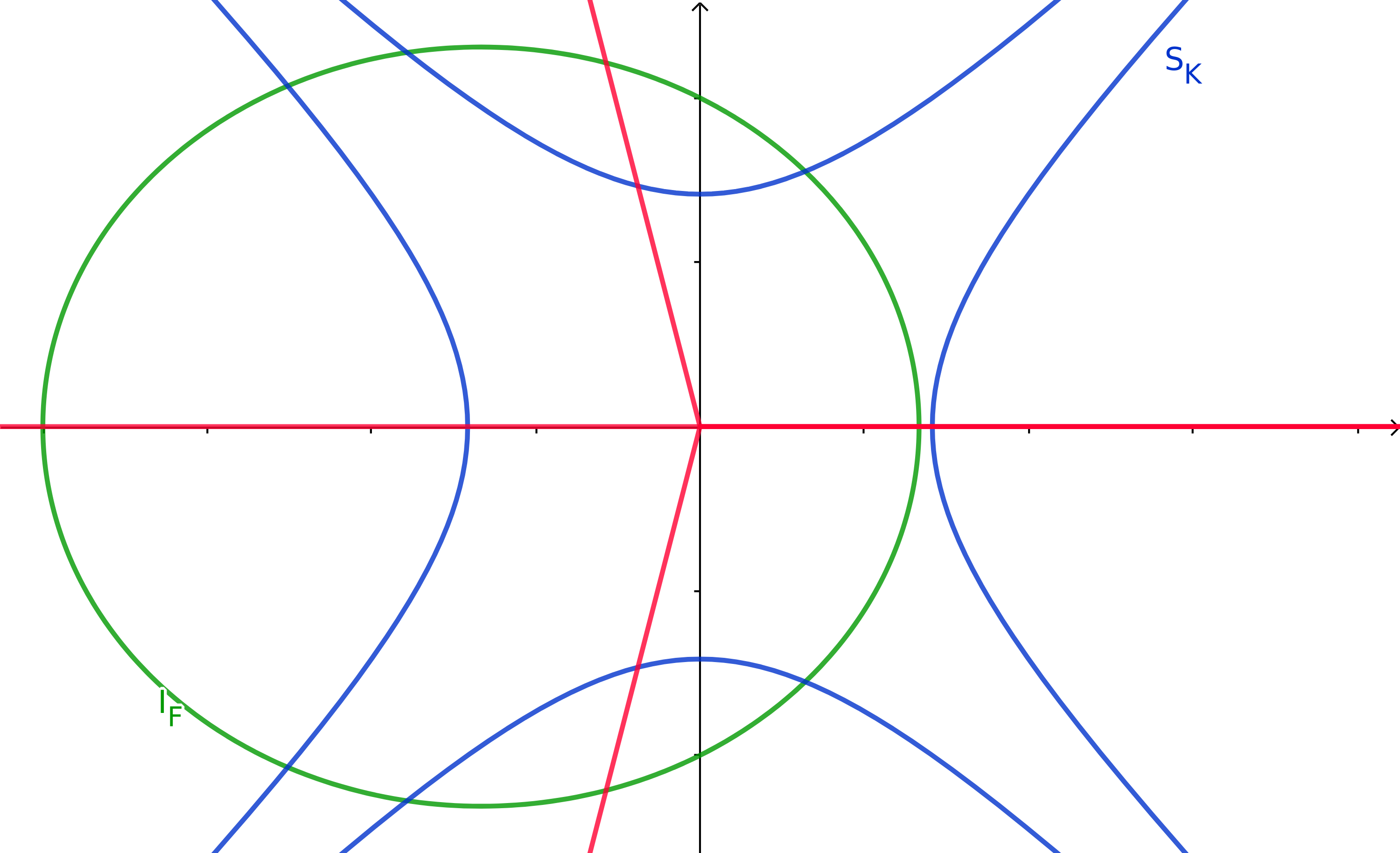}

\section{The existence}
We state the main theorem about the existence of homogeneous geodesics,
depending on the signature of the Killing form.
We recall that the first part of theorem was stated in \cite{DuS}, however, the proof was given
with the assumption that the Killing form is semidefinite.
We present here the complete proof which shows that the same idea for finding two geodesic vectors
works on each component $K=1$ and $K=-1$ of the unit (pseudo-)sphere of the Killing form.
Let us remark that the present result is optimal. The examples given above illustrate that this
result cannot be improved in general. The examples of solvable Lie groups with invariant
Finsler metrics which admit just two homogeneous geodesics were given in \cite{DuSr2}.

As we will not express vectors in components now, it is more convenient to change notation and
use lowercase letters for vectors.
\begin{theorem}
Let $(M,F)$ be a homogeneous Finsler manifold.
There exist at least two homogeneous geodesics through arbitrary point $p\in M$.
Let $K$ be the Killing form of a transitive isometry group $G$ of $M$
and let $\fg=\fh+\fm$ be a reductive dedomposition.
If the restriction of $K$ to $\fm$ is indefinite,
then there exist at least four homogeneous geodesics through arbitrary point $p\in M$.
\end{theorem}
{\it Proof.}
Let $G$ be a transitive isometry group of $M$ and let $H$ be the isotropy group of a fixed point $p\in M$.
We express $M$ as the homogeneous space $M=G/H$.
Let $K$ be the Killing form of $G$ and let $\rad(K)$ be the null space of $K$.
We choose $\fm=\fh^\perp$ with respect to $K$.
The decomposition in ${\mathrm{Ad}}(H)$-invariant and the Finsler metric induces the invariant Minkowski norm
and its fundamental tensor on $\fm$. We shall denote these again by $F$ and $g$.
The Killing form $K$ negatively definite on $\fh$, because $H$ is compact.
Hence, $\rad(K)\subseteq\fm$. We shall distinguish the two cases:

Case 1) $\rad(K)=\fm$: We choose a hyperplane $W\subset\fm$ such that $[\fm,\fm]\subset W$.
There exist two vectors $n_1,n_2\in\fm$ which are normal to $W$, which means
\begin{eqnarray}
\nonumber
g_{n_i} ( n_i, w ) & = & 0 \qquad \forall w\in W.
\end{eqnarray}
Vectors $n_1, n_2$ are on the different sides of the hyperplane $W$ and they
are in general not opposite to each other, unless $F$ is reversible.
See \cite{S} or \cite{DuS} for details of the construction.
We shall now write $n$ for any of the two vectors $n_1,n_2$.
For an arbitrary fixed vector $w\in W$,
the function $F^2(n+tw)$ attains its minimum at $t=0$ and hence,
using the standard formula, we obtain
\begin{eqnarray}
\nonumber
0 = \frac{1}{2} \frac{d}{dt} F^2(n+tw) \big |_{t=0} = g_n(n,w), \qquad \forall w\in W,
\end{eqnarray}
which is the desired property. In particular, it is satisfied for any $w\in [\fm,\fm]\subset W$.
We obtain immediately, using Lemma \ref{golema2}, that $n_1$ and $n_2$ are geodesic vectors.

Case 2) $\rad(K)\subsetneq\fm$: We start with the construction and notation as in \cite{YH}.
We shall investigate the function
\begin{eqnarray}
\nonumber
f(z) & = & \frac{K(z,z)}{F^2(z)}.
\end{eqnarray}
This function is positively homogeneous and it is reasonable to restrict
the definition domain to the indicatrix
\begin{eqnarray}
\nonumber
I_F=\{z\in\fm; F(z)=1\}.
\end{eqnarray}
Since the group $H$ is compact and $\rad(K)$ is an ${\mathrm{Ad}}(H)$-invariant subspace,
there exist an ${\mathrm{Ad}}(H)$-invariant $K$-orthogonal complement $W$ of $\rad(K)$ in $\fm$.
Each vector $z\in\fm$ can be uniquely decomposed as $z=z_1+z_2$, where $z_1\in\rad(K)$ and $z_2\in W$.
Denote $k=\dim(\rad(K))$ and let
\begin{eqnarray}
\nonumber
D_k=\{z_1\in\rad(K), F(z_1)<1\}
\end{eqnarray}
be the open unit disc in $\rad(K)$.
For each fixed $z_1\in D_k$, consider the submanifold
\begin{eqnarray}
\nonumber
S_{z_1} = \{ z_2\in W, F(z_1+z_2)= 1 \}\subset W,
\end{eqnarray}
which has the topology of the sphere ${\mathbb{S}}^{m-1}$, where $m={\mathrm{dim}}(W)={\mathrm{dim}}(M)-k$.
From now on, $z_1+z_2$ means $z_1\in D_k$, $z_2\in S_{z_1}$ and $z_1+z_2\in I_F$.
Each sphere $S_{z_1}\subset W$ is split by the nullcone of $K|_W$ into open submanifolds
$S_{z_1}^+ = \{ z_2\in S_{z_1} ; K(z_2,z_2)>0 \}$ and
$S_{z_1}^- = \{ z_2\in S_{z_1} ; K(z_2,z_2)<0 \}$.
Both boundaries $\partial \bar S_{z_1}^+$ and $\partial \bar S_{z_1}^-$
of closures $\bar S_{z_1}^+$ and $\bar S_{z_1}^-$
are the intersection of $S_{z_1}$ with the nullcone of $K|_W$.
It is easily seen that the function $f(z)$ defined above is positive on $S_{z_1}^+$ and negative on $S_{z_1}^-$.
Limits of $f(z_1+z_2)$ for $z_2\in S_{z_1}^+$, $z_2\rightarrow \partial \bar S_{z_1}^+$ and for
$z_2\in S_{z_1}^-$, $z_2\rightarrow \partial \bar S_{z_1}^-$ and are all zero.
For the later use, we define already now the distinguished open submanifolds of $I_F$, namely
\begin{eqnarray}
\nonumber
S^+ = \{ z_1+z_2\in I_F ; K(z_2,z_2)>0 \} = \cup_{z_1\in D_k} (z_1+S_{z_1}^+), \cr
S^- = \{ z_1+z_2\in I_F ; K(z_2,z_2)<0 \} = \cup_{z_1\in D_k} (z_1+S_{z_1}^-).
\end{eqnarray}
We further denote by $K_1$ and $K_{-1}$ the standard (pseudo-)spheres with respect to the Killing form $K|_W$,
namely $K_1= \{ w\in W ; K(w,w)=1 \}$ and $K_{-1}= \{ w\in W ; K(w,w)=-1 \}$.
It is easy to observe that, for fixed $z_1$, for each vector $z_2\in S_{z_1}^+$,
there is a positive real number $p$ such that $pz_2\in K_1$.
In the same way, for each vector $z_2\in S_{z_1}^-$ there is a positive real number $p$ such that $pz_2\in K_{-1}$.
This correspondence gives homeomorphisms $\psi^+_{z_1}$ of each $S_{z_1}^+$ with $K_1$
and homeomorphisms $\psi^-_{z_1}$ of each $S_{z_1}^-$ with $K_{-1}$.

Let the signature of $K$ be $(p,q,k)$, where $p$ is the number of positive signs in the diagonal form of $K$,
$q$ is the number of negative signs  and $p+q=m$.
It is well known (see for example \cite{BON}) that the topology of $K_1$, or $K_{-1}$, respectively,
(and hence also topology of each $S_{z_1}^+$, or each $S_{z_1}^-$, respectively) is the topology of
${\mathbb{S}}^{p-1}\times{\mathbb{R}}^{q}$, or ${\mathbb{R}}^{p}\times {\mathbb{S}}^{q-1}$, respectively.
In the special case $q=0$, or $p=0$, respectively, it reduces to the topology of the sphere ${\mathbb{S}}^{m-1}$.
In the special case $p=1$, or $q=1$, respectively, it reduces to the topology of
the two copies of ${\mathbb{R}}^{q}$, or the two copies of ${\mathbb{R}}^{p}$, respectively.
We continue with the general case $p>1, q>1$.
We investigate the manifold $S^+$, which is homeomorphic to
$D_k\times {\mathbb{S}}^{p-1} \times {\mathbb{R}}^{q}$
and the function $f(z)$ defined above is positive on it.
We have observed that, on each $S_{z_1}^+$, $\lim_{z_2\rightarrow \partial \bar S_{z_1}^+} f(z_1+z_2) = 0$.
It is also easy to see that for $z_1\in D_k$, it holds $\lim_{z_1\rightarrow \partial \bar D_k} f(z_1+z_2) =0$.
Obviously, the function $f(z)$ reaches its maximum on $S^+$ for some vector $y_1\in S^+$.
Now we are going to show that there exist a vector $y_2\in S^+$,
where the function $f(z)$ on $S^+$ reaches the saddle point.

We identify $W$ with ${\mathbb{R}}^{p} \times {\mathbb{R}}^{q}$ and we fix the homeomorphism
$\phi\colon {\mathbb{S}}^{p-1}\times{\mathbb{R}}^{q} \rightarrow K_1\subset {\mathbb{R}}^{p} \times {\mathbb{R}}^{q}$
by the formula $(s,x)\mapsto (\sqrt{1+|x|^2}\cdot s,x)$.
We define, for fixed $z_1\in {\mathrm{rad}}(K)$ and fixed $x\in {\mathbb{R}}^{q}$, the submanifolds
$C_{z_1,x}$ of $S_{z_1}^+$ as $((\psi^+_{z_1})^{-1}\circ\phi)({\mathbb{S}}^{p-1}\times x)$.
Each $C_{z_1,x}$ is homeomorphic to the sphere ${\mathbb{S}}^{p-1}$.
For fixed $z_1$ and $x$, and with compact definition domain $C_{z_1,x}$, the function $f(z_1+z_2)$ restricted to $C_{z_1,x}$
attains its minimum $\varepsilon(z_1,x)>0$ at some $z_1+\bar z_2(z_1,x)\in C_{z_1,x}$.
For each $z_1\in D_k$ and $x\in {\mathbb{R}}^{q}$, we choose one such $\bar z_2$ and consider the mapping
$\varphi\colon D_k\times{\mathbb{R}}^q \rightarrow I_F$, $(z_1,x)\mapsto z_1+\bar z_2(z_1,x)$.
The function $f(\varphi(z_1,x))=\varepsilon(z_1,x)$ is smooth on $D_k\times{\mathbb{R}}^q$
and it attains its maximum $\lambda_2$ at $(\bar z_1,\bar x)$.
Here $(\bar z_1,\bar x)$ can be choosen and the map $\varphi$ can be defined in a way that there is
a neighbourhood $U\subset D_k\times{\mathbb{R}}^q$ of $(\bar z_1,\bar x)$
such that the mapping $\varphi\big |_U$ is smooth.
We put $y_2 = \varphi(\bar z_1,\bar x)\in S^+$.
The definition of $y_2$ does not depend on the identification of $W$
with ${\mathbb{R}}^{p} \times {\mathbb{R}}^{q}$.
From the construction, it follows the existence of a basis
$B=\{ u_1,\dots,u_{p-1},v_1,\dots,v_k,w_1,\dots,w_q\}$
of $T_{y_2}S^+$ with the following property:
Vectors $u_i$ form a basis of the tangent space of $C_{\bar z_1,\bar x}$ at $y_2=\varphi(\bar z_1,\bar x)$,
vectors $v_i$ are images in the tangent mapping to $\varphi$ of tangent vectors to $D_k$ at $(z_1,x)$ and
vectors $w_i$ are images in the tangent mapping to $\varphi$ of tangent vectors to ${\mathbb{R}}^q$ at $(z_1,x)$.
The function $f(z)$ attains its local minimum along any curve $\gamma(t)$ in $S^+$
with $\gamma(0)=y_2$ and whose tangent vector at $t=0$ is any of the vectors $u_i$.
And the function $f(z)$ attains its local maximum along any curve $\gamma(t)$ in $S^+$
with $\gamma(0)=y_2$ and whose tangent vector at $t=0$ is any of the vectors $v_i$ or $w_i$.

It remains to show that $y_1$ and $y_2$ are geodesic vectors.
As to $y_1$, the function
\begin{eqnarray}
\nonumber
\tilde f(z) = K(z,z) - \lambda_1 F^2(z)
\end{eqnarray}
attains its minimum $0$ at $y_1$.
For any fixed $w\in \fm$, the function $\hat f(t)=\tilde f(y_1+tw)$ attains its minimum $0$ at $t=0$
and hence $\hat f'(0)=0$.
It follows that
\begin{eqnarray}
\label{eq0}
K(y_1,w) & = & \lambda_1 \cdot g_{y_1}(y_1,w), \qquad \forall w\in \fm.
\end{eqnarray}
As to $y_2$, the function
\begin{eqnarray}
\nonumber
\tilde f(z) = K(z,z) - \lambda_2 F^2(z)
\end{eqnarray}
attains value $0$ at $y_2$.
For any vector $u_i$ defined above, the function $\hat f(t)=\tilde f(y_2+tu_i)$ attains its maximum $0$ at $t=0$
and hence $\hat f'(0)=0$.
For any of the vectors $v_i$ or $w_i$ defined above, the function $\hat f(t)=\tilde f(y_2+tu_i)$
attains its minimum $0$ at $t=0$ and hence also $\hat f'(0)=0$.
It follows that
\begin{eqnarray}
\label{eq1}
K(y_2,w) & = & \lambda_2 \cdot g_{y_2}(y_2,w)
\end{eqnarray}
holds for any vector $w$ from the above basis $B$ of $T_{y_2}S^+$.
It is obvious that this equality holds also for $w=y_2$ and consequently formula (\ref{eq1})
holds for any $w\in\fm$.
Formulas (\ref{eq0}) and (\ref{eq1}) lead to formula
\begin{eqnarray}
\nonumber
g_{y_i}({y_i},[y_i,z]_\fm) = \frac{1}{\lambda_i} K(y_i,[y_i,z]_\fm) = \frac{1}{\lambda_i} K([y_i,y_i],z) = 0,
 \qquad \forall z\in\fm,\, i=1,2,
\end{eqnarray}
which shows that $y_1$ and $y_2$ are geodesic vectors.

Finally, similar construction with $S^-$ and with the function $-f(z)$ leads to geodesic vectors $y_3, y_4\in S^-$.
If $p=1$, the manifold $S^+$ has two connected components and vectors $y_1, y_2$ are chosen
as vectors where the function $f(z)$ reaches its maximum on each of these components.
The saddle point may not exist in this situation, as the examples above illustrate.
If $q=1$, the manifold $S^-$ has two connected components and vectors $y_3, y_4$ are chosen in similar way.
If $p=0$ or $q=0$, then either $S^+$ or $S^-$ is trivial and the procedure leads only to two geodesic vectors.
$\hfill\square$

\section*{Acknowledgements}
The research was supported by the grant IGS 8210-017/2020 of the Internal Grant Agency
of Institute of Technology and Business in \v Cesk\'e Bud\v ejovice.

\end{document}